\documentclass[10pt]{amsart}

\numberwithin{equation}{section}

\newcommand{\ben}{\begin{enumerate}}
\newcommand{\een}{\end{enumerate}}

\newcommand{\bea}{\begin{eqnarray}}
\newcommand{\ba}{\begin{array}}
\newcommand{\bean}{\begin{eqnarray*}}
\newcommand{\ea}{\end{array}}
\newcommand{\eea}{\end{eqnarray}}
\newcommand{\eean}{\end{eqnarray*}}
\newcommand{\beq}{\begin{equation}}
\newcommand{\eeq}{\end{equation}}
\newcommand{\bthm}{\begin{thm}}
\newcommand{\ethm}{\end{thm}}
\newcommand{\blem}{\begin{lem}}
\newcommand{\elem}{\end{lem}}
\newcommand{\bprop}{\begin{prop}}
\newcommand{\eprop}{\end{prop}}
\newcommand{\bcor}{\begin{cor}}
\newcommand{\ecor}{\end{cor}}
\newcommand{\bdfn}{\begin{dfn}}
\newcommand{\edfn}{\end{dfn}}
\newcommand{\brem}{\begin{rem}}
\newcommand{\erem}{\end{rem}}
\newcommand{\bpf}{\begin{proof}}
\newcommand{\epf}{\end{proof}}
\newcommand{\bfact}{\begin{fact}}
\newcommand{\efact}{\end{fact}}
\newcommand{\bobs}{\begin{obs}}
\newcommand{\eobs}{\end{obs}}

\alph{enumii} \roman{enumiii}

\newtheorem{thm}{Theorem}[section]
\newtheorem{prop}[thm]{Proposition}
\newtheorem{lem}[thm]{Lemma}

\newtheorem{cor}[thm]{Corollary}
\newtheorem{dfn}[thm]{Definition}
\newtheorem{rem}[thm]{Remark}
\newtheorem{fact}[thm]{Fact}
\newtheorem{obs}[thm]{Observation}

\def\cH{\mathcal H}             \def\cF{\mathcal F}       
                   \def\cP{{\mathcal P}}

\def\C{{\mathbb C}}                      \def\oc{{\hat \C}}

\newcommand{\amsc}{{\mathbb C}}

\newcommand{\cbar}{\hat{{\mathbb C}} }

\def\a{\alpha}                \def\b{\beta}             \def\d{\delta}
               \def\e{\varepsilon}          
\def\g{\gamma}                \def\Ga{\Gamma}           \def\l{\lambda}
              \def\om{\omega}           
               \def\sg{\sigma}
                          
\def\ka{\kappa}

\newcommand{\ep}{\varepsilon}
\newcommand{\ph}{\varphi}

\def\1{1\!\!1}

\def\and{\text{ and }}

        \def\diam{\text{\rm {diam}}}
\def\dist{\text{{\rm dist}}}

     \def\HD{\text{{\rm HD}}}

\def\bi{\bigcap}              \def\bu{\bigcup}
\def\({\bigl(}                \def\){\bigr)}
\def\lt{\left}                \def\rt{\right}

\def\ld{\ldots}                        \def\^{\tilde}

\def\es{\emptyset}            \def\sms{\setminus}
\def\sbt{\subset}             \def\spt{\supset}
\def\gek{\succeq}             \def\lek{\preceq}

\def\comp{\asymp}
           \def\downto{\searrow}
\def\sp{\medskip}                     
\def\ov{\overline}            

\def\ni{\noindent}

\def\om{\omega}

\def\supp{\text{{\rm supp}}}

\newcommand{\jul}{J_f}
\newcommand{\fat}{{\mathcal F}_f}
\newcommand{\sing}{sing(f^{-1})}

\newcommand{\cri}{{\mathcal C}_f}
\newcommand{\post}{{\mathcal P}_f}
\newcommand{\av}{{\mathcal A}_f}

\newcommand{\pft}{{\mathcal{L}}_t}
\newcommand{\pfh}{{\mathcal{L}}_h}
\newcommand{\pfs}{{\mathcal{L}}_s}
\newcommand{\npfs}{\tilde{\mathcal{L}}_s}

\begin{document}
%***********************************************

\title[]
{ \bf\large {\Large E}rgodic properties of sub-hyperbolic functions with polynomial Schwarzian derivative}
\date{\today}
% % author information %
\author[\sc Volker MAYER]{\sc Volker MAYER}
%\\[0.5cm]\rm \large Preliminary Version}
\address{Volker Mayer, Universit\'e de Lille I, UFR de Math\'ematiques,
UMR 8524 du CNRS,
59655 Villeneuve d'Ascq Cedex, France}
\email{volker.mayer\@math.univ-lille1.fr\newline \hspace*{0.3cm} Web:
math.univ-lille1.fr/$\sim$mayer}
\author[\sc Mariusz URBA\'NSKI]{\sc Mariusz URBA\'NSKI}
\address{Mariusz Urba\'nski, Department of Mathematics,
 University of North Texas, Denton, TX 76203-1430, USA}
\email{urbanski\@unt.edu\newline \hspace*{0.3cm} Web:
www.math.unt.edu/$\sim$urbanski}
%
% dedication %
%\dedicatory{}
%
% AMS information %
\thanks{Research of the second author supported in part by the
NSF Grant DMS 0400481.}

\begin{abstract}
The ergodic theory and geometry of the Julia set of meromorphic functions
on the complex plane
with polynomial Schwarzian derivative is investigated under the condition
that the forward trajectory of asymptotic values in the Julia set is
bounded and the map $f$ restricted to its closure is expanding, the
property refered to as sub-expanding.
We first show the existence, uniqueness, conservativity and ergodicity
 of a conformal measure $m$ with minimal exponent $h$; furthermore, we show
weak metrical exactness of this measure. Then we prove
 the existence of a $\sg$--finite invariant measure $\mu$ absolutely
 continuous with respect to $m$. Our main result states that $\mu$ is finite
 if and only if the order $\rho$ of the function $f$ satisfies the
 condition $h>3\frac{\rho}{\rho +1}$. When finite, this measure is shown to be
metrically exact. We also establish a version of Bowen's
 formula showing that the exponent $h$ equals the Hausdorff dimension of the Julia set
 of $f$.

\end{abstract}

\maketitle

%***********************************************
\section{Introduction}

The study of the ergodic theory and geometry of the Julia set
of transcendental meromorphic functions appears to be a delicate task
due to the infinite degree of these functions. For example, even the
existence of conformal measures, on which the whole theory relies and
which is by now completely standard in the realm of
 rational functions or Kleinian groups, is not known in general.
Employing Nevanlinna's theory and a convenient change of the Riemannian
metric we provided a complete treatise for a very general class of
hyperbolic meromorphic functions in the papers \cite{mu1} and \cite{mu2}.
In the present paper we relax the hyperbolicity assumption and allow
the Julia set to contain singularities. Clearly one can adopt the arguments
developed in the theory of rational iteration to deal with certain type
of critical points. More challenging is to analyze the contribution of logarithmic
singularities and, as we will see, this gives quite surprising results.
The class of meromorphic functions with polynomial Schwarzian derivatives fit best
to such a project since they have logarithmic singularities but they do not
have critical points. We therefore restrict our considerations to this class
of functions which, in particular, contains the tangent family;
definitions and other examples are given in Section~2.

In the context of ergodic theory and fractal
geometry, meromorphic and entire functions with logarithmic
singularities have been investigated in \cite{sk1}, \cite{sk2}, \cite{uzdnh}
(see also \cite{ku} for a more complete historical outline and list of references)
and, more recently, in \cite{ks}. In
\cite{sk1} and \cite{sk2} these singularities landed at poles an, in \cite{uzdnh}, they
were escaping to infinity extremely (like the trajectory of zero under
the exponential function) fast. In both of these cases the forward
trajectory of images of logarithmic singularities experienced a
large expansion neutralizing the contracting effect of singularities themselves.
Assuming that a meromorphic map is subhyperbolic, the postsingular set is bounded,
the Julia set is an entire sphere, and the reference conformal measure is the
Lebesgue measure, the paper \cite{ks} addressed the role of logarithmic singularities
(algebraic singularities were also allowed).

In the present paper we consider subexpanding
 meromorphic functions $f:\C\to\oc$ with polynomial
Schwarzian derivative. By subexpanding we understand
that the postsingular set $\post$ is bounded, and that the map $f$
restricted to $\post$ is expanding.
Employing the full power of Nevanlinna theory we first prove the existence
of an atomless conformal measure via the Patterson-Sullivan construction.
This measure is proved to be weakly metrically exact, which implies its
ergodicity and conservativity. We then show the following result in which the existence of
the $\sg$--finite measure $\mu$ is obtained by employing M. Martens general method.

\

\bthm \label{thm main}
Let $f$ be a subhyperbolic meromorphic function $f$ of polynomial Schwarzian derivative
and let $m$ be the $h$--conformal measure of $f$ obtained via the Patterson--Sullivan construction.
Then there exists a $\sg$--finite invariant measure $\mu$
absolutely continuous with respect to $m$. Moreover, the measure
$$
\mu \;\text{ is finite}\quad \text{ if and only if} \qquad h>3\frac{\rho}{\rho +1}
$$
 where $\rho = \rho (f)$ is the order of the function $f$. If $\mu$ is finite, then the
dynamical systems ($f,\mu$) it generates is metrically exact and, in consequence, its Rokhlin's
natural extension is K-mixing.
\ethm

\

Notice that $3\frac{\rho}{\rho +1}\geq 2$ if and only if the order
$\rho \geq 2$. Consequently the measure $\mu$ is most often infinite. However,
in the case of the
tangent family, which is just one specific example among others, this invariant
measure can be finite. Curiously,
finiteness of the invariant measure
for the strictly preperiodic function $z\mapsto 2\pi i e^z$ is not
known as yet.

Let us mention that we do not assume that the Julia set is the entire
sphere nor that the conformal measure is the Lebesgue measure. In fact we do not
assume that any conformal measure exists at all. But in the special situation
when the Julia set is the entire sphere (in which case the spherical Lebesque measure
is automatically a conformal measure) and if in addition $h=2>3\frac{\rho}{\rho +1}$,
i.e. if the order of the function $\rho <2$, then the existence of a probability
invariant measure absolutely continuous with respect to the Lebesgue measure
follows also from \cite{ks}.
Indeed, in that situation our necessary and sufficient condition $\rho <2$
conincides with the sufficient condition (Z3) from the paper by Kotus and \'Swiatek.
Concerning the reciprocal statement, \cite{ks} simply provides a counterexample.

The most involved part of the proof of Theorem~\ref{thm main} is to
show finiteness. In the case the measure $\mu$
is finite, the dynamical system it generates is shown to be $K$--mixing which, 
in particular, implies mixing of all orders.

We also investigate the Hausdorff dimension of the Julia set and show that this 
dimension coincides with
$h$, the exponent of the conformal measure $m$. Notice that this holds
despite that the $h$--dimensional Hausdorff measure
is shown to vanish on the Julia set.

\

\section{The class of functions and definitions}

\subsection{Definitions}

 The reader may consult, for example,
\cite{nev'}, \cite{nev} or \cite{hille} for a detailed exposition on
meromorphic functions and \cite{b} for their dynamical aspects. We
collect here the properties of interest for our concerns. The \it Julia
set \rm of a meromorphic function $f:\amsc\to \cbar$ is denoted by
$\jul$ and the \it Fatou set \rm by $\fat$. Note that, in contrast to
\cite{mu1, mu2}, we include here $\infty \in \jul$
since we are dealing with spherical geometry.
However, $O^{-}(\infty )$ is a very special subset of the Julia set.

Let $\av $ be the set of \it asymptotic values. \rm Note that the functions we consider
do not have critical values. Therefore
$\av$ coincides with the so called set of \it singular values \rm $\sing$.
The \it post-singular set \rm $\post$ is the closure (in the sphere) of the set
$\bigcup_{n>0} f^n(\av )$.\\

Concerning the singularities of a meromorphic function $f$, we dispose of Iversen's classification
(see e.g. \cite{b}): let $a\in \sing $ and, for every $r>0$,
$U_r$ be a component of $f^{-1}(D(a,r))$ in such a way that
$r_1<r_2$ implies $U_{r_1}\subset U_{r_2}$. Then there are two
possibilities:
\begin{itemize}
  \item[a)] $\bigcap _{r>0} U_r=\{c\}$ consists of one point, or
  \item[b)] $\bigcap _{r>0} U_r=\emptyset$.
\end{itemize}
In the latter case we say that our choice $r\mapsto U_r$ defines a
\it transcendental singularity of $f^{-1}$ over $a$. \rm Such a
singularity is called \it logarithmic \rm if the restriction
$f:U_r\to D(a,r)\setminus \{a\}$ is a universal cover for some
$r>0$. If this is the case, then the component $U_r$ is called \it
logarithmic tract. \rm For the functions we consider all the
transcendental singularities are logarithmic.

In case a), the point $c$ can be regular or it is a critical
point $c\in \cri$.\\

 We will always denote by
$$d\sg (z)=\frac{|dz|}{1+|z|^2}$$ the spherical metric and by
$$|f'(z)|_\sg =|f'(z)| \frac{1+|z|^2}{1+|f(z)|^2}$$
the derivative of $f$ with respect to the spherical metric.
The following direct consequence of Koebe's distortion theorem
will be used.

\blem\label{bdp}
Let $f:\C\to\oc$ be a meromorphic function and suppose that
$D(w,2\d )\subset \oc \setminus \post$. Then, for every $n\geq 1$,
$z\in f^{-n}(w)$ and all $x,y\in D(w,\d )$ we have that
$$ K^{-1} \leq \frac{|(f_z^{-n})'(y)|_\sg}{|(f_z^{-n})'(x)|_\sg} \leq K$$
for some universal constant $K\geq 1$.
\elem

Here and in the rest of the paper $f_z^{-n}$ signifies the inverse branch
of $f^n$ defined near $f^n (z)$ mapping back $f^n(z)$ to $z$. An other convention
will be that $D(z,r)$ stands for the
spherical metric centered at $z$ and of radius $r$.
To indicate a spherical $r$--neighborhood of a set $X$ we write $B(X,r)$.

\subsection{Meromorphic functions with polynomial Schwarzian derivative}

We consider meromorphic functions $f:\C\to\oc$ for which the
Schwarzian derivative
\beq \label{2.4} S(f) = \left( \frac{f''}{f'} \right)' -\frac{1}{2} \left(
  \frac{f''}{f'} \right)^2= 2P\eeq
is a polynomial and for which the set of asymptotic values $\av$
does not contain infinity.
Nevanlinna \cite{nev3} established that meromorphic functions with
polynomial Schwarzian derivative
are exactly the functions that have only finitely many asymptotical
values and no critical values.
Moreover, if such a function has a pole, then it is of order
one. Consequently the maps of this class
are locally injective. We also mention that any solution of
(\ref{2.4}) is of order $\rho = p/2$,
where $p=deg (P) +2$, and it is of normal type of its order
(cf. \cite{hille2}).

 Standard examples are furnished by the
tangent family $f(z)=\l \tan (z)$ for which $S(f)$ is constant. By M\"obius invariance of $S(f)$, functions like
$$\frac{e^z}{\l e^z +e^{-z}} \quad and \quad \frac{\l e^z}{e^z -e^{-z}} $$
have also constant Schwarzian derivative. Examples for which $S(f)$ is a polynomial are
\begin{equation} \label{2.0.1}
f(z) =\frac{a\,Ai(z) +b\,Bi(z)}{c\,Ai(z)+d\,Bi(z)} \qquad with \quad
ad-bc\neq 0
\end{equation}
where $Ai$ and $Bi$ are the Airy functions of the first and second kind.
These a linear independent solutions of $g''-zg=0$ and, in general, if
$g_1,g_2$ are linear independent
solutions of
\begin{equation}\label{2.3.2}
g'' + Pg =0 \; ,
\end{equation}
then $f=\frac{g_1}{g_2}$ is a solution of the Schwarzian equation (\ref{2.4}).
Conversely, every solution of (\ref{2.4}) can be written locally as a
quotient of two linear independent solutions
of the linear differential equation (\ref{2.3.2}). The asymptotic properties of these solutions
are well known due to work of Hille (\cite{hille3}, see also \cite{hille2}). They give a
precise description of the function $f$ near infinity. We now collect
the facts that are important for our needs (more details and references are for example in \cite{mu2}).

First of
all, there are $p$ critical directions $\theta_1,...,\theta_p$ which
are given by
$$arg \, c +p\theta =0\, (mod \, 2\pi )$$
where $c$ is the leading coefficient of $P(z)=cz^{p-2} +...$. In a sector
$$ S_j =\Big\{ |arg z - \theta _j | < \frac{2\pi}{p}-\d \; ; \; |z| >R\Big\}\, ,$$
 $R>0$ is sufficiently large and $\d > 0$, the equation (\ref{2.3.2}) has two linear independent solutions
\begin{equation}\label{2.3.5}
\begin{array}{l}
g_1 (z) = P(z) ^{-\frac{1}{4}} exp \big(  iZ +o(1)  \big) \;\quad and \\
g_2 (z)= P(z) ^{-\frac{1}{4}} exp \big(  -iZ +o(1) \big)
\end{array}
\end{equation}
where
\beq\label{2.3.8}Z=\int_{2Re^{i\theta_j}}^z P(t)^\frac{1}{2} \, dt = \frac{2}{p}c^\frac{1}{2} z^\frac{p}{2}\big(
1+o(1)\big) \quad for \quad z\to \infty \;\; in \;\; S_j.\eeq
If $f$ is a meromorphic solution of the Schwarzian equation (\ref{2.4}), then there are $a,b,c,d\in \C$
with $ad-bc\neq 0$ such that
\begin{equation}\label{2.3.6}
f(z) = \frac{ag_1(z)+bg_2(z)}{cg_1(z) +dg_2(z)} \quad , \;\; z\in S_j .
\end{equation}
Observe that $f(z) \to a/c$ if $z\to \infty$ on any ray in
$S_j\cap \{arg \, z <\theta _j \}$ and that $f(z)\to b/d$
if $z\to \infty$ on any ray in $S_j\cap \{arg \, z >\theta _j \}$.
The asymptotic values of $f$ are given by all the
$a/b$, $c/d$ corresponding to all the sectors $S_j$, $j=1,...,p$.

With this precise description of the asymptotic behavior of $f$,
one can show (\cite{mu2}) that
\beq\label{2.3.9}
|f'(z)|\asymp |z|^{\rho-1}\big|\a +\b f(z) +\g f^2(z)\big| \quad for \;\; z\in S_j
\eeq
where $\a =-ab/\d$, $\b =(ad+bc)/\d$, $\g =-cd/\d$ and $\d =ad-bc $. Notice that
$\g\neq 0$ if all the asymptotic values
of $f$ are finite.\\

\subsection{Sub-hyperbolic functions}

\bdfn \label{2.0}
The function $f$ is called boundedly non-recurrent if $\infty \not\in\av\cup\post$, if
 $\av \cap \post \cap \jul=\emptyset$ and if every asymptotic value that belongs
to the Fatou set is in an attracting component. If in addition, $f$ is
expanding on $\post$, the map $f$ is called sub-expanding.
\edfn

\

\ni Notice that this definition implies that all the asymptotic values
of the function $f$ are finite and that the post-singular set is
bounded and nowhere dense in the Julia set. From now on we fix a
number $T>0$ such with the following properties.
\begin{itemize}
\item[(T$_1$)] $4T<|a_1-a_2| $ for all distinct $a_1,a_2\in \av$.
\item[(T$_2$)] $B(\post ,4T) \cap \av =\emptyset$.
\item[(T$_3$)] $f^{-1}(\av\cup\post ) \cap \big( B(\av\cup\post ,
  4T)\setminus (\av \cup \post )\big)
  =\emptyset$.\\
\end{itemize}

\

\ni To every asymptotic value $a\in\av$ there correspond (finitely many ) logarithmic tracts $U_a$.
In the following such a tract $U_a$ will always be a component of $D(a,T)$ and we may
suppose that $U_a \cap B(\av \cup\post , 4T)=\emptyset$.\\

\sp\ni Notice that being boundedly non-recurrent implies sub-expanding
whenever the Julia set is equal to $\oc$ (\cite{gks}). We do always
assume this property. From now on we also require $T>0$ to be so
small  that $|(f^p)'|_\sg >2$ on $B(\post,T)$ for some $p\geq 1$ and
that  there are open neighborhoods
$\Omega_1, \Omega_0$ of $\post$ such that $\overline{\Omega_1} \subset \Omega_0\subset B(\post , T)$ and
$g=f^p_{|\Omega_1}:\Omega_1 \to \Omega _0$ is a proper mapping. Denote
\beq\label{2.5} \Omega_n=g^{-n}(\Omega _0) \quad
and\quad \Gamma _n = \Omega_n \setminus \Omega_{n+1}.\eeq

\

\ni Using the facts that repelling periodic points are dense in $\jul$ and that
$\jul$ contains poles, one can easily prove the following.

\

\bobs\label{o1v54}(topological exactness of $f$) For every non-empty open set $U$
intersecting $\jul$, there exists $n\ge 0$ such that $f^n(U)\spt\oc\sms\av$. In
particular, for every $r>0$ there exists $q_r\ge 0$ such that $f^{q_r}(D(z,r))
\spt\oc\sms\av$ for all $z\in\jul$.
\eobs

\

\ni Since $\post$ is a closed forward-invariant set and the map $f|_{\post}$ is expanding,
following inverse trajectory of a point near $\post$, one can prove the following.

\

\bobs\label{o2v54} (repeller)
The set $\post$ is a repeller for $f$, precisely, assuming
$T>0$ to be small enough, we have
$$
\bi_{n=0}^\infty f^{-n}(B(\post,2T))=\post.
$$
\eobs

\subsection{First observations and transfer operator}
If one choses the right metric space $(\C, d\sg)$, then the ergodic theory
of meromorphic functions
can be well developped. This has been done in great generality and in the
hyperbolic case in \cite{mu1, mu2}.
For the functions we consider here the right geometry is simply the spherical
one (which result from (\ref{2.3.9}). Indeed,
the functions satisfy the balanced growth condition of \cite{mu1} with $\a_1=\rho -1$ and
with $\a_2=2$ the later meaning that one has to work with the spherical metric).

\

\blem \label{2.1}
Let $f:\C\to\oc$ be of polynomial Schwarzian derivative with $\infty \not\in \av$.
Then, if $z$ belongs to a logarithmic tract $U_{2T}\subset f^{-1}(D(a,2T))$ over an asymptotic
value $a\in \av$, we have that
$$
|f'(z)|_\sg\asymp (1+|z|^{\rho + 1})|f(z) -a|
$$
and otherwise
$$
|f'(z)|_\sg\asymp 1+|z|^{\rho + 1}
$$
where $\rho=\rho(f)<\infty$ is the order of $f$.
\elem

\bpf
Follows from asymptotic description of $f$ near infinity, in particular (\ref{2.3.9}),
together with the fact that $f$ has only simple poles.
\epf

\

Let us consider the transfer operator with respect to the spherical geometry.
\beq \label{2.2}
\pft \ph (w)= \sum_{z\in f^{-1}(w) } |f'(z)|_\sg ^{-t} \ph(z) \quad \ph \in C(\jul).
\eeq
 It follows from Lemma \ref{2.1} that
\beq \label{2.2.1} \pft \1 (w) \preceq \max\{1, dist(w,\av )^{-t}\}
\sum_{z\in f^{-1}(w) }(1+|z|^{\rho + 1}) ^{-t}\eeq
for every $w\in \oc \setminus \av  $.
This last sum is very well known in the theory of meromorphic functions and
a Theorem of Borel \cite{nev} together with the divergence
property of $f$ established in Theorem 3.2 of \cite{mu2} implies that
\beq \label{2.3} \pft \1 (w)<\infty \quad \text{if and only if } \quad t >\frac{\rho}{\rho +1}.
\eeq
We need the following additional properties.

\

\bprop \label{2.0.2} For every
$t>\frac{\rho}{\rho +1}$, there exists a constant $M_t$ such that
$$\Sigma (t,w)=\sum_{z\in f^{-1}(w) }(1+|z|^{\rho + 1}) ^{-t}
\leq M_t \quad \text{for every } \;\; w\in \oc .$$
\eprop

\

\ni The proof of this result uses parts of \cite{mu1,mu2} and relies
heavily on Nevanlinna theory. Good references for this are
\cite{nev'} or \cite{cy}. Let us simply recall that $n_f(r,a)$
stands for the number of $a$--points of modulus at most $t$, that
the integrated counting number $N_f(r,a)$ is defined by
$dN_f(r,a)=n_f(r,a)/r$ and that $T_f(r)$ denotes the characteristic
function of $f$.

\bpf
 Fix $\ep =1$ and let $A>0$ a constant that will be
precised
later. We may suppose that the origin is not a pole of $f$.

{\it Case 1.} $w\not\in D(f(0), \ep )$. Then we have that
$$\Sigma(t,w) \preceq  C_A + \sum_{\begin{array}{c}
                  f(z)=w \\
                  |z|>A
                \end{array}}(1+|z|^{\rho + 1} )^{-t}
                \preceq C_A + \sum_{\begin{array}{c}
                  f(z)=w \\
                  |z|>A
                \end{array}}|z|^{-u}
                $$
with $C_A = \sup_{w\in \oc} n_f(w,A)<\infty$ and with $u=(\rho + 1) t$.
Since $f$ is of finite order $\rho $ we can make the following two
integrations by part:
$$
 \sum_{\begin{array}{c}
                  f(z)=w \\
                  |z|>A
                \end{array}}|z|^{-u}
                = -\frac{n_f(A,w)}{A^u} - u \frac{N_f(A,w)}{A^{u+1}} +
                u^2 \int_A ^\infty \frac{N_f(s,w)}{s^{u+1}}
                \leq u^2 \int_A ^\infty \frac{N_f(s,w)}{s^{u+1}}.
$$
 The First Main Theorm of Nevanlinna theorem (see
\cite[Corollary 4.2]{mu1}) gives that
$$ N_f(r,w) \leq T_f(r) -\log[f(0),w],$$
where $[a,b]$ denotes the chordal distance on the Riemann sphere
(with in particular $[a,b]\leq 1$ for all $a,b\in\oc$). Since in
this first case $w\not\in D(f(0), \ep )$ there is $\Theta <\infty$
such that
$$ N_f(r,w) \leq T_f(r) + \Theta \quad for \;\; every \;\; w\not\in  D(f(0), \ep ).$$
Therefore
$$\sum_{\begin{array}{c}
                  f(z)=w \\
                  |z|>A
                \end{array}}(1+|z|^{\rho +1 } )^{-t} \preceq u^2
                \int_A^\infty \frac{T_f(s)+\Theta}{s{u+1}} ds
                =\tilde{M}_{u} <\infty\quad for \;\; every \;\; w\not\in D(f(0), \ep ).$$
All in all, there exists $M^{(1)}_u <\infty$ such that
$$ \Sigma(t,w) \leq M^{(1)}_u \quad for \;\; every \;\; w\not\in  D(f(0), \ep ).$$

{\it Case 2.} $w\in  D(f(0), \ep )$. We are let to find a
uniform bound for
$$\sum_{\begin{array}{c}
                  f(z)=w \\
                  |z|>A
                \end{array}}|z|^{-u} \quad , \;\; w\in \jul \cap D(f(0), \ep ).$$

Let $v\in \C$ be a point that is not a pole of $f$ and such that
$|f(-v)-f(0)|>2\ep$. Set $A= 3|v|$ and define the meromorphic
function $g(\xi )=f(\xi-v)+v$. If $\xi =z+v$, then $f(z)=w$ is
equivalent to $g(\xi)=w+v$.

Notice that $g(0)=f(-v)+v$. If we set $a=w+v$, then
\beq \label{2.0.5}
|a-g(0)| = |w-f(-v)| \geq |f(-v)-f(0)| -|f(0)-w| > \ep.\eeq
 On the
other hand, if $|\xi-v|\geq A$, then $|\xi| \geq A-|v| \geq 2|v|$ and
$\frac{1}{|\xi-v|} \leq \frac{2}{|\xi |}$.
It follows that
$$ \sum_{\begin{array}{c}
                  f(z)=w \\
                  |z|>A
                \end{array}}|z|^{-u}=\sum_{\begin{array}{c}
                  g(\xi )=a \\
                  |\xi-v|>A
                \end{array}}|\xi -v|^{-u}\leq
                 \sum_{\begin{array}{c}
                  g(\xi)=a \\
                  |\xi|>2|v|
                \end{array}}\left(\frac{2}{|\xi|}\right)^{u}.$$
In the same way as before we can now use again the First Main
Theorem of Nevanlinna theory, this time applied to the function $g$.
Remember that by (\ref{2.0.5}) we have $a\not\in D(g(0), \ep )$
whenever $w=a-v\in D(f(0), \ep )$. Therefore,
$$\sum_{\begin{array}{c}
                  f(z)=w \\
                  |z|>A
                \end{array}}|z|^{-u} \leq 2^u
                \sum_{\begin{array}{c}
                  g(\xi)=a \\
                  |\xi|>2|v|
                \end{array}}|\xi|^{-u}\leq \tilde{\tilde{M}}_u$$
                for every $a= w+v$, $w\in  D(f(0), \ep )$.
It follows that there is $M^{(2)}_u <\infty$ such that
$$ \Sigma (t,w) \leq M^{(2)}_u \quad for \;\; every \;\; w\in  D(f(0), \ep ).$$
All in all we showed that there is $M_t <\infty$ such that $\| \Sigma(t,.)
\| _\infty \leq M_t$. \epf

\section{Conformal measures}
A probability measure $m$ is a $t-$conformal measure for the meromorphic function
$f:\C\to \oc$ if $\frac{dm\circ f}{dm}=|f'|_\sg^t$ or, equivalently, if for every
mesurable set $E\subset \C$ for which the restriction $f_{|E}$ is injective we have
$$ m(f(E)) = \int _E |f'|_\sg^t \, dm .$$
Notice that if $\jul=\oc$, which is the case when all the asymptotic values are strictly preperiodic,
the spherical Lebesgue measure is a $2$--conformal measure.
Therefore we restrict in the following subsection to functions with non-empty Fatou set.

\subsection{Existence of conformal measures and the pressure function}
Conformal measures are usually obtained via the standard Patterson-Sullivan
method (see for example \cite{mcm} which contains a nice description of this procedure).
 For meromorphic functions however one must very carefully check what
is going on at infinity because at this point the function is not defined.
We ignore this for a moment and work on the compact set $\jul \subset \oc$.
Let us first consider the Poincar\'e series of $f$ at $\infty$
\beq \label{3.1}
\cP(t)=\cP(t,\infty)
= \sum_{n>0}\sum_{z\in f^{-n}(\infty)} |(f^n)'(z)|_\sg^{-t}
=\sum_{n>0} \pft ^n \1 (\infty)
\eeq
which is well defined for $t>\frac{\rho}{\rho+1}$.
Define
\beq \label{3.2}
h =h_f=  \inf \left\{ t> \frac{\rho}{\rho+1}\;  ; \quad \cP (t) <\infty \right\}.
\eeq
We have the following.

\

\blem \label{3.3}
For a sub-hyperbolic meromorphic function $f$ of polynomial Schwarzian derivative
with $\fat \neq \emptyset$ we have that
$$ \frac{\rho}{\rho+1}< h_f \leq 2.$$
\elem

\bpf
The function $f$ cannot have Baker nor rotation domains (see \cite{b}). Therefore one can
find a disc $D\subset \fat$ arbitrarily close to $\infty$ such that
all the inverse images $f^{-n}(D)$, $n\geq 0$,
are disjoint. It follows then from Koebe's distortion theorem that
$$
\cP(2)=\sum_{n>0}\sum_{z\in f^{-n}(\infty)} |(f^n)'(z)|_\sg^{-2}
\asymp \sum_{n>0} meas (f^{-n} (D)) <\infty,
$$
where $meas$ stands for the spherical Lebesgue measure.
The fact that $h>\frac{\rho}{\rho+1}$
is an immediate consequence of \cite{my} together with the divergence property of $f$ given
in Theorem 3.2 of \cite{mu2}.
\epf

\

\ni The Patterson-Sullivan method as described in \cite{mcm} applies now
and furnishes a $h$--conformal measure say $m$. Using this measure we can make
the following important improvement of the above estimate for $h_f$.

\

\bprop\label{3.4}
For a sub-hyperbolic meromorphic function $f$ of polynomial Schwarzian derivative
with $\fat \neq \emptyset$
we have that
$$
2\frac{\rho}{\rho+1}< h_f .
$$
\eprop

\bpf
Let $a\in \av$ and let $a'=f(a)$. The function $f$ is expanding on $\post$: there is
$p>0$ such that $|(f^p)'|_\sg >2$ on $B(\post , T)$. Let again $g=f^p_{\Omega_1}$ (see \ref{2.5}).
 Notice that $g(D(g^n(a'),T))\supset D(g^{n+1}(a'),T)$.
Denote
\beq\label{3.5}
A_n=D(g^n(a'),T) \setminus g_{g^n(a')}^{-1} \Big(D(g^{n+1}(a'),T)\Big).
\eeq
Increasing $p$ if necessary we have from the fact that $\post$ is compact and nowhere dense in $\jul$ that
$$\inf _{n\geq 0} m(A_n) \geq c>0 .$$
If $V_n =f_a^{-1}\circ g_{a'}^{-n} (A_n)$ then
\beq\label{3.6}
m(V_n) \asymp \Big( |f'(a)|_\sg |(g^n)'(a')|_\sg \Big)^{-h} m(A_n) \asymp diam(V_n)^h \quad , \quad n\geq 1.
\eeq
The preimages of $V_n$ to a logarithmic tract $U$ over $a$ can be labeled by
$U_{n,k} =f_k^{-1} (A_n)$. Let $z_{n,k}\in U_{n,k}$ be any point. Then
\beq\label{3.4a}
\aligned
m(U_{n,k})
&\asymp |f'(z_{n,k})|_\sg^{-h}diam(V_n)^h \\
&\asymp |z_{n,k}|^{-(\rho +1)h} |f(z_{n,k})-a|^{-h} diam(V_n)^h
\asymp |z_{n,k}|^{-(\rho+1)h},
\endaligned
\eeq
where the relation
\beq\label{3.4b}
|z_{n,k}|\asymp(n^2+k^2)^{{1\over 2\rho}}
\eeq
follows from an elementary calculation based on (\ref{2.3.8}) and (\ref{2.3.6}). Hence,
\beq\label{3.7}
1\geq m(U) =\sum_{n,k} m(U_{n,k})
\asymp \sum_{n,k} |z_{n,k}|^{-(\rho+1) h}
\asymp \sum_{n,k} \left(n^2+k^2 \right)^{-\frac{\rho+1}{2\rho } h}
\eeq
The assertion of the Lemma follows now since this last sum is convergent if and only if
$\frac{\rho+1}{2\rho } h >1$.
\epf

\

\ni As a first application of this estimate on $h_f$ we can now prove that the measure $m$ does
not charge infinity.

\

\blem \label{3.8}
For the $h$--conformal measure $m$ obtained from the Patterson-Sullivan construction we have $m(\{\infty \})=0$.
\elem

\bpf
The measure $m$ is obtained in the following way. If the Poincar\'e series $\cP(t)$ diverges for $t=h$,
then $m$ is a limit of measures of the form
$$\nu_s = \frac{1}{\cP(s)} \sum_n \sum_{z\in f^{-n}(\infty )} |(f^n)'(z)|_\sg ^{-s} \delta_z
= \frac{1}{\cP(s)} \sum_n (\pfs ^n)^*\delta_\infty .$$
If ever $\cP(h)<\infty$ then one must add artificially a divergence behavior.
A simple way of doing this is to follow the exposition of \cite{mcm} and to replace the exponent $s=h+\delta$ by $s'=h-\delta$
  in a finite number (depending on $s$) of terms in this expression. To take into account these modifications we write
  $$\nu_s = \frac{1}{\tilde{\cP}(s)} \sum_n(\npfs ^n)^*\delta_\infty .$$
 Notice that $\delta\to 0$. We can therefore
  suppose that $s'=h-\delta \geq \tau > 2\frac{\rho}{\rho +1}$. For every $E\subset \C$
  for which $f_{|E}$ is injective one has
  \beq\label{3.9}\int_E \min\big\{|f'|_\sg ^s , |f'|_\sg ^{s'} \big\}d\nu_s \leq \nu_s (f(E)) \leq
  \int_E \max\big\{|f'|_\sg ^s , |f'|_\sg ^{s'} \big\}d\nu_s \eeq
  (see \cite{mcm} for details). By classical arguments (that the reader can find in \cite{du1}) it follows that
  \beq\label{3.10} \nu_s (E) \leq \int _\C \sum_{z\in f^{-1}(w)}\max\big\{|f'(z)|_\sg ^{-s} , |f'(z)|_\sg ^{-s'} \big\} \1_E (z) \, d\nu_s(w) \eeq
  for every mesurable set $E\subset \C$.

 By definition $\nu_s(\{\infty\})=0$. We have to show that the sequence $(\nu_s)_s$ is tight at infinity
 which means that for every $\e>0$ there is $R>0$ such that $\nu_s(W_R)<\e$ for every $s>h$ where
   $W_R =\{|z|>R\}$. This set can be written as
   $$W_R = \bigcup_{a\in\av} (U_a \cap W_R ) \cup \tilde{W}_R$$
   where the union is taken over all the (finitely many) logarithmic tracts $U_a \subset f^{-1}(D(a,T)$ over
   the asymptotic values $a\in\av$.

   Tightness on $\tilde{W}_R$ can be shown like in \cite{mu1}. The key point is the following.
   If $z\in \tilde{W}_R$ then $w=f(z)\not\in B(\av , T)$. Therefore, if $\g>0$ is any small number such that
   $(\rho+1)\tau  -\g >\rho $, then it follows from (\ref{3.10}) and from Proposition \ref{2.3} that
   \begin{eqnarray*}
   % \nonumber
     \nu_s(\tilde{W}_R ) &\preceq &
\int _\C \sum_{z\in f^{-1}(w)}\1 _{\tilde{W}_R} (z) |z|^{-(\rho +1)\tau} d\nu_s(w)\\
      &\leq& \frac{1}{R^\gamma} \int _\C \sum_{z\in f^{-1}(w)} (1+|z|)^{\g -(\rho+1)\tau}d\nu_s(w)
\leq \frac{M}{R^\g}.
   \end{eqnarray*}
Let us now consider what happens on a logarithmic tract $U=U_a \subset f^{-1}(D(a,T))$ over $a\in \av$.
With the notations of the proof of Proposition \ref{3.2} and with the same arguments,
one has $\nu_s(V_n) \preceq diam(V_n)^\tau$
and
$$ \nu_s(U_{n,k}) \preceq |f'(z_{n,k})|_\sg ^{-\tau} diam (V_n)^\tau\asymp |z_{n,k}|^{-(\rho +1)\tau}
 \asymp \left(n^2+k^2\right)^{-\frac{\rho+1}{2\rho}\tau}.$$
 Notice that these estimates do not depend on $s>h$ and imply that
 $$
\nu_s(U) = \sum_{n,k} \nu_s(U_{n,k})\preceq \sum_{n,k}\left(n^2+k^2\right)^{-\frac{\rho+1}{2\rho}\tau}
\leq c<\infty
$$
 since $\tau>2\frac{\rho}{\rho +1}$. Therefore $\nu_s(U\cap W_R) \to 0$ as $R\to \infty$ uniformly in $s>h$.
\epf

\subsection{Additional properties}

Recall the definition of the annuli $\Gamma_n$ are given in (\ref{2.5}). We start with the
following.

\

\blem \label{3.10}
There exists $0<\g<1$ such that $m(\Gamma _n )\preceq \g ^n $ for every $n \geq 0$.
\elem

\bpf
Let $\post \subset \bigcup_{j=1}^N D_j$ where the discs $D_j =D(x_j,2T)$, $x_j\in \post$,
built a Besicovitch covering of $\post$. We may suppose that $\Omega_0\subset\bigcup_{j=1}^N D_j$.
Fix $q\geq 1$ such that for every $j=1,...,N$
$$ m(\Gamma _q \cap D_j ) \leq \eta \, m(\Gamma _0 \cap D_j)$$
with $\eta>0$ some small number to be determined later on. Remember that $g=f^p_{|\Omega _1}$.
Clearly all the inverse branches of $g^n$ are well defined and of bounded distortion on
every disc $D_j$. Let us denote these by $g_*^{-n}$.
With this notation we can calculate, for every $n\geq 1$, that
\begin{eqnarray*}
  m(g^{-n} (D_j\cap \Gamma _q)) &=& \sum_*  m(g^{-n}_* (D_j\cap \Gamma _q))
=\sum_*  \frac{ m(g^{-n}_* (D_j\cap \Gamma _q))}{ m(g^{-n}_*
(D_j\cap \Gamma _0))} m(g^{-n}_* (D_j\cap \Gamma _0))\\
   &\preceq& \sum_*  \frac{ m(D_j\cap \Gamma _q)}{ m(D_j\cap \Gamma _0)}
m(g^{-n}_* (D_j\cap \Gamma _0))
\leq \eta \sum_*m(g^{-n}_* (D_j\cap \Gamma _0))\\
   &=& \eta \, m(g^{-n} (D_j\cap \Gamma _0)).
\end{eqnarray*}
Summing over $j$ and using the Besicovitch property of the covering we get that
$$ m(\Gamma _{q+n}) \leq C \eta \, m(\Gamma _n) \quad for \;\; every \; n\geq 0.$$
The assertion follows provided $\eta$ has been chosen such that $C\eta <1/2$.
\epf

\

\ni In the rest of this section we denote $\nu$
any $h$--conformal measure (and keep the letter $m$
for the conformal measure that has been constructed above).
Note that for any Borel probability measure $\nu$ on a compact metric space
$(X,\rho)$,
$$
M_\nu(r):=\inf\{\nu(B(x,r):x\in\supp(\nu)\}>0
$$
for every $r>0$.
Let us also prove the following.

\

\blem\label{l3v54}
For any $h$--conformal measure $\nu$ we have $\nu(\post)=0$.
\elem

\

\bpf Recall that one condition imposen on $T$ was that for every $z\in\post$ and
every $n\ge 0$, there exists a holomorphic inverse branch $f_z^{-n}:D(f^n(z),2T)
\to\oc$ of $f^n$ sending $f^n(z)$ to $z$. It then follows from the bounded
distortion property (Lemma ~\ref{bdp}) that
\beq\label{1v54}
\aligned
\nu(D(z,K^{-1}T|(f^n)'(z)|_\sg^{-1}))
&\le \nu\(f_z^{-n}(D(f^n(z),T))\)\\
& \preceq |(f^n)'(z)|_\sg^{-h}\nu(D(f^n(z),T))
\le |(f^n)'(z)|_\sg^{-h}.
\endaligned
\eeq
Since $\post$ is a nowhere dense subset of $\jul$, there exists $\g>0$ such that
for every $y\in\post$ there exists $\hat y\in\jul$ such that
$$
D_y:=D(\hat y,\g)\sbt D(y,K^{-2}T)\sms\post.
$$
Then
\beq\label{1v54.1}
f_z^{-n}(D_{f^n(z)})
\sbt f_z^{-n}(D(f^n(z),K^{-2}T))\sms\post
\sbt D(z,K^{-1}T|(f^n)'(z)|_\sg^{-1})\sms\post
\eeq
and
$$
\nu\(f_z^{-n}(D_{f^n(z)})\)\succeq |(f^n)'(z)|_\sg^{-h}\nu (D_{f^n(z)})\ge M_\nu (\g)|(f^n)'(z)|_\sg^{-h}.
$$
Combining this, (\ref{1v54}), (\ref{1v54.1}), and noting that
$\supp(\nu )=\jul$, we get that
$$
{\nu \(D(z,K^{-1}T|(f^n)'(z)|_\sg^{-1})\sms\post\)\over
\nu \(D(z,K^{-1}T|(f^n)'(z)|_\sg^{-1})\)}
\succeq M_\nu (\g) \quad \text{for every } \quad n\geq 1.
$$
Therefore,
$$
\limsup_{r\to 0}{\nu (D(z,r)\sms\post)\over \nu (D(z,r))}
\succeq  M_\nu (\g)>0.
$$
So, $z$ is not a Lebesgue density point of $\nu $, and therefore $\nu (\post)=0$.
\epf

\subsection{Metric exactness, conservativity and ergodicity}
Suppose that $(X,\cF,\nu)$ is a probability space and $T:X\to X$ is a measurable map
such that $T(A)\in\cF$ whenever $A\in\cF$. The map $T:X\to X$ is said to be \it weakly
metrically exact \rm  provided that $\ov{\lim}_{n\to\infty}\nu(T^n(A))=1$ whenever $A\in\cF$
and $\nu(A)>0$. A straightforward observation concerning weak metrical
exactness is this.

\

\bobs\label{o1v54.1}
If a measurable transformation $T:X\to X$ of a probability space $(X,\cF,\nu)$ is
weakly metrically exact, then it is ergodic and conservative.
\eobs

\

\ni In the context of invariant measures there is the following, more involved fact, also
indicating a dynamical significance of weak
metrical exactness (see e.g. \cite{cfs, pu}).

\

\bfact\label{t1v54.3}
A measure-preserving transformation $T:X\to X$ of a probability space $(X,\cF,\mu)$ is
weakly metrically exact if and only if it is exact, which means that $\lim_{n\to\infty}
\mu(T^n(A))=1$ whenever $A\in\cF$ and $\mu(A)>0$, or equivalently, the $\sg$-algebra
$\bi_{n\ge 0}T^{-n}(\cF)$ consists of sets of measure $0$ and $1$ only. Then the
Rokhlin's natural extension $(\^T,\^X,\^\mu)$ of $(T,X,\mu)$ is K-mixing.
\efact

\

\ni The main result of this subsection is this.

\

\bthm\label{t1v55}
$m$ is a unique $h$-conformal measure. The dynamical system $f:\jul\to\jul$ is
weakly metrically exact with respect to $m$. In particular it is ergodic and conservative.
\ethm

\bpf
Let
$$
\post^*=\{z\in\jul:\dist_\sg(z,\av\cup\post)>2T\}.
$$
By Observation \ref{o2v54},
\beq\label{2v55}
\jul^*=\{z\in\jul\sms O^-(\infty):\om(z)\cap\post^*\ne\es\}
=\jul\sms\bu_{n=0}^\infty f^{-n}(\post\cup\{\infty\}).
\eeq
Take $z\in\jul^*$. Then there exists a strictly increasing sequence
$(n_j=n_j(z))_{j=1}^\infty$ of positive integers such that
$$
f^{n_j}(z)\in\post^*\sms\{\infty\}
$$
for all $j\ge 1$. Then for every $j\ge 1$, there exists a meromorphic
inverse branch $f_z^{-n_j}:D(f^{n_j}(z),2T)\to\oc$ of $f^{n_j}$ sending
$f^{n_j}(z)$ to $z$. It then follows from Lemma~\ref{bdp} (bounded
distortion property) that for every $h$-conformal measure $\nu$ on $\jul$,
\beq\label{3v55}
\aligned
\nu(D(z,K^{-1}T|(f^{n_j})'(z)|_\sg^{-1}))
&\le \nu\(f_z^{-n_j}(D(f^{n_j}(z),T))\)\\
& \preceq |(f^{n_j})'(z)|_\sg^{-h}\nu(D(f^{n_j}(z),T))
\le |(f^{n_j})'(z)|_\sg^{-h}.
\endaligned
\eeq
Put
$$
r_j(z)=(4K)^{-1}T|(f^{n_j})'(z)|_\sg^{-1}.
$$
The above formula rewrites then as follows.
\beq\label{1v55}
\nu(D(z,4r_j(z)))\preceq r_j^h(z).
\eeq
It also follows from Lemma~\ref{bdp} that
\beq\label{1v57}
\aligned
\nu(D(z,r_j(z)))
&\ge \nu\(f_z^{-n_j}(D(f^{n_j}(z),(4K^2)^{-1}T))\) \\
&\succeq |(f^{n_j})'(z)|_\sg^{-h}\nu(D(f^{n_j}(z),4K^2)^{-1}T)) \\
&\ge M_\nu((4K^2)^{-1}T)|(f^{n_j})'(z)|_\sg^{-h} \\
&\asymp r_j^h(z).
\endaligned
\eeq
Now fix $E$, an arbitrary Borel set contained in $\jul^*$. Fix also $\e>0$. Since
the measure $m$ is regular, for every $z\in E$ there exists $j(z)\ge 1$ such that,
with $r(z)=r_{j(z)}(z)$, we will have
\beq\label{2v57}
m\lt(\bu_{z\in E}D(z,r(z))\rt)\le m(E)+\e.
\eeq
By the $(4r)$-covering theorem there exists now a countable set $\hat E\sbt E$
such that the balls $\{D(z,r(z))\}_{z\in\hat E}$ are mutually disjoint and
$$
\bu_{z\in\hat E}D(z,4r(z))\spt \bu_{z\in E}D(z,r(z))\spt E.
$$
Hence, using (\ref{1v55}), (\ref{1v57}) (with $\nu$ replaced by $m$) and
(\ref{2v57}), we get
$$
\aligned
\nu(E)
&\le \sum_{z\in\hat E}\nu(D(z,4r(z)))
 \le (4K^2/T)^h\sum_{z\in\hat E}r^h(z) \\
&\le K^{2h}M_m((4K^2)^{-1}T)\sum_{z\in\hat E}m(D(z,r(z))) \\
&\asymp m\lt(\bu_{z\in\hat E}D(z,r(z)\rt) \le m(E)+\e.
\endaligned
$$
Thus, letting $\e\downto 0$, we get $\nu(E)\preceq m(E)$.
Hence, $\nu_{|\jul^*}$ is absolutely continuous with respect to $m_{|\jul^*}$.
Exchanging the roles of
$\nu$ and $m$ we get $m_{|\jul^*}\lek\nu_{|\jul^*}$, and finally that $\nu_{|\jul^*}$
is equivalent to $m_{|\jul^*}$.
Since, in view of Lemma~\ref{l3v54},
$$m\lt(\bu_{n=0}^\infty f^{-n}(\post )\rt)
=\nu\lt(\bu_{n=0}^\infty f^{-n}(\post )\rt)=0,$$ we thus conclude that $\nu$ and $ m$
are equivalent on $\jul \setminus O^-(\infty )$, $O^-(\infty ) = \bu_{n=0}^\infty f^{-n} (\{\infty\})$.
Finally, if $\nu \left( O^-(\infty )\right) >0$ then $\nu^* =\nu _{| O^-(\infty )}$
would be a conformal measure without mass on $\jul \setminus O^-(\infty )$.
But then we would have a contradiction since we have just seen that $m$ and $\nu^*$ are equivalent on
$\jul \setminus O^-(\infty )$. Therefore $\nu \left( O^-(\infty )\right) =0$ and both measures
are equivalent on the whole Julia set.

Passing to the proof of weak metrical exactness of $f$ with respect to the measure
$m$, suppose that $E\sbt\jul$ and
\beq\label{1v59}
\limsup_{n\to\infty}\sup\{m(f^n(E)\cap D(y,K^{-2}T))/m(D(y,K^{-2}T)):y\in\post^*\}=1.
\eeq
We shall show that
\beq\label{2v59}
\limsup_{n\to\infty}m(f^n(E))=1.
\eeq
In virtue of Observation~\ref{o1v54} there exists $q\ge 0$ such that
$$
f^q(D(y,K^{-2}T))\spt\oc\sms \av
$$
for all $y\in\jul$. Clearly, by conformality of $m$, for every $\e>0$ there then
exists $\d>0$ such that if $y\in\jul$, $G\sbt D(y,K^{-2}T)$, and $m(G)/m(D(y,K^{-2}T))\ge 1-\d$,
then $m(f^q(G))\ge 1-\e$. Combining this with (\ref{1v59}) yields (\ref{2v59}).
In order to the weak metrical exactness of $m$, suppose by contrapositive that
$E\sbt \jul$ and $\limsup_{n\to\infty}m(f^n(E))<1$. By (\ref{1v59}) and (\ref{2v59}),
this implies that
$$
2\ka:=\liminf_{n\to\infty}\inf\lt\{m(D(y,K^{-2}T\sms f^n(E))/m(D(y,K^{-2}T)):y\in\post^*\rt\}>0.
$$
So, for all $n\ge 1$ large enough, say $n\ge p$,
$$
\inf\lt\{m(D(y,K^{-2}T\sms f^n(E))/m(D(y,K^{-2}T)):y\in\post^*\rt\}\ge\ka>0.
$$
Fix $z\in E\cap \jul^*$. We shall show that $z$ is not a Lebesgue density point
for the measure $m$. Let $n_j=n_j(z)\ge p$, $j\ge 1$, have the same meaning as in the
first part of the proof. Then
\beq\label{3v59}
f_z^{-n_j}\(D(f^{n_j}(z),K^{-2}T)\sms f^{n_j}(E)\)
\sbt D(z,K^{-1}T|(f^{n_j})'(z)|_\sg^{-1})\sms E
\eeq
and
$$
\aligned
m\(f_z^{-n_j}\(D(f^{n_j}(z) &,K^{-2}T)\sms f^{n_j}(E)\)\)\ge \\
&\ge K^{-h}|(f^{n_j})'(z)|_\sg^{-h})m\(D(f^{n_j}(z),K^{-2}T)\sms f^{n_j}(E)\) \\
&\ge \ka K^{-h}|(f^{n_j})'(z)|_\sg^{-h})m\(D(f^{n_j}(z),K^{-2}T)\) \\
&\ge \ka K^{-h}M_m(K^{-2}T)|(f^{n_j})'(z)|_\sg^{-h}).
\endaligned
$$
Combining this along with (\ref{3v59}) and (\ref{3v55}), we get that
$$
{m\(D(z,K^{-1}T|(f^{n_j})'(z)|_\sg^{-1})\sms E\)
 \over m\(D(z,K^{-1}T|(f^{n_j})'(z)|_\sg^{-1})\)}
\ge \ka K^{-2h}M_m(K^{-2}T).
$$
Therefore,
$$
\lim_{r\to 0}{m(D(z,r)\sms E)\over m(D(z,r))}\ge \ka K^{-2h}M_m(K^{-2}T)>0.
$$
So, $z$ is not a Lebesgue density point for $m$. Thus $m(E\cap\jul^*)=0$. Since
$m(\jul^*)=1$ (see (Lemma~\ref{l3v54} and (\ref{2v55}), we finally get that $m(E)
=0$. The weak metrical exactness of $f$ with respect to $m$ is established.
Ergodicity and conservativity follow from Observation~\ref{o1v54.1}. Since $\nu$
(introduced in the first part of the proof) is equivalent to $m$, the equality
$\nu=m$ follows from ergodicity of $m$. We are done. \epf

\section{Invariant Measures}
We now consider $f:\C\to\oc$ a sub-hyperbolic meromorphic function $f$
of polynomial Schwarzian derivative
and investigate invariant measures equivalent to the conformal measure $m$
obtained in the previous section. In particular we show in the course of this section
Theorem \ref{thm main}.

\subsection{Existence of $\sg$--finite invariant measures.}
Since we already established conservativity of the conformal measure $m$ we can
use the method of M. Martens \cite{mm} (see also \cite{ku} for a description of this method)
 in order to obtain the following.

\

 \bprop\label{4.1}
Let $f$ be a sub-hyperbolic meromorphic function $f$ of polynomial Schwarzian derivative
and let $m$ be the conservative $h$--conformal measure of $f$ with $m(\post )=m(\{\infty \})=0$.
Then there exists $\mu$ a $\sg$--finite invariant measure
absolutely continuous with respect to $m$.
 \eprop

\bpf
Using a Whitney decomposition of $\C\setminus (\av \cup \post )$ it is easy to construct a countable
partition $\{ A_n  ; \, n\geq 0\}$ of $X=\jul \setminus (\{\infty \} \cup \av \cup \post)$ such that
 for every  $n,m\geq 0$ there is $k\geq 0$ such that
$$
m(f^{-k}(A_m) \cap A_n ) >0 .
$$
Since $m$ has no mass on $\jul\setminus X$ and since $m$ is conservative
M. Martens result \cite{mm} applies and gives the $\sg$--finite invariant measure
absolutely continuous with respect to $m$. Notice that for every Borel set $A\subset X$ we have that
\beq \label{4.2}
\mu (A) = \lim_{n\to \infty } \frac{\sum_{k=0}^n m(f^{-k}(A))}{\sum_{k=0}^n m(f^{-k}(A_0))}.
\eeq
For the choice of the set $A_0$ there is much freedom. We will use in particular that $A_0\subset X$ is
such that all the inverse of the iterates of $f$ are well defined and have bounded distortion.
\epf

\

\ni Let $\Delta = \oc \setminus B(\av \cup \post , T)$.

\

\blem \label{4.3}
There is $K>1$ such that $1/K \leq \ph <K$ on $\Delta$ where $\ph =d\mu / dm$.
\elem

\bpf
Let $z\in \Delta $.
From the expression (\ref{4.2}) follows that
$$ \ph (z) =\lim_{r\to 0}\frac{\mu (D(z,r))}{m (D(z,r))} \asymp\lim_{r\to 0}\frac{1}{m (D(z,r))}
 \lim_{n\to \infty } \frac{\sum_{k=0}^n \pfh^k\1 (z) m(D(z,r))}{\sum_{k=0}^n \pfh^k\1(z_0)m(A_0)}$$
 where $z_0\in A_0$ is any point.
 Now, if $z_1, z_2\in\Delta $ are any two points, then they can be joined by a chain of at most $N=N(\Delta )$
spherical discs of radius $T$. On each of these discs all the inverse branches of every
iterate of $f$ is well defined and have distortion bounded by some universal constant.
Therefore
$$\pfh^k\1 (z) \asymp \pfh^k\1 (z_0) \quad for \;\; every \; k\geq 0.$$
The Lemma is proven.
\epf

\

\ni This simple observation on the density $h$ has several important applications
starting with the following.

\

\bprop \label{4.4}
$\mu (B(\av , T)) <\infty$.
\eprop

\bpf
It suffices to show that $\mu (D(a,T))<\infty$, $a\in \av$. The measure $\mu$ being invariant,
$\mu (D(a,T)) =\mu (f^{-1}(D(a,T))) $.
By the choice of the constant $T>0$ (see (T$_3$)), $f^{-1}(D(a,T))\cap B(\av \cup \post , T) =\emptyset$.
It therefore follows from Lemma \ref{4.3} that
$$\mu (D(a,T)) =\mu (f^{-1}(D(a,T))) \asymp m(f^{-1}(D(a,T)) <\infty.$$
\epf

\subsection{When is the $\sg$--finite invariant measure finite?}

\ni To our big surprise it turns out that finiteness of the invariant measure $\mu$
does depend on the order of the function.

\

\bthm\label{4.5}
Let $f$ be a sub-hyperbolic meromorphic function of polynomial Schwarzian derivative
and let $m$ be the (unique) $h_f$--conformal measure of $f$. Then there is a finite $f$--invariant
measure $\mu$ absolutely continuous with respect to $m$ if and only if $h>3 \frac{\rho}{\rho +1}$.
\ethm

\

\ni Consequently the invariant measure $\mu$ is finite in the particular case of the tangent family
and also for the examples of (\ref{2.0.1}) that involve the Airy functions.
Notice that $3 \frac{\rho}{\rho +1}\geq 2$ as soon as the order $\rho =deg(P)/2 +1 \geq 2$.

\

\bcor \label{4.6}
A sub-hyperbolic solution $f$ of the polynomial Schwarzian equation $S(f)=2P$ has
finite invariant measure absolutely continuous with respect to the $h_f$--conformal measure
if and only if $deg(P)= 0$ or $deg(P)=1$.
\ecor

\

\ni We now prove Theorem \ref{4.5} in several steps and use again the notations given in (\ref{2.5}).
 Let us consider in the following $f$ a sub-hyperbolic
meromorphic function of polynomial Schwarzian derivative
and let $m$ be again the $h_f$--conformal measure of $f$.

\

\blem\label{4.7}
If $h_f \leq 3 \frac{\rho}{\rho +1}$ then there is no finite invariant measure absolutely
continuous with respect to $m$.
\elem

\bpf
Suppose to the contrary that such a finite invariant measure $\mu$ exists.
Remember that $g=f^p_{|\Omega _1} :\Omega_1\to \Omega _0$. Let $a\in \av$, set $a'' =f^p(a)$
and $D''=D(a'', T)$. By invariance of $\mu$ we have that
$$ \mu (\Omega _n ) = \mu(f^{-p} (\Omega _n ))
\geq \mu (f^{-p}_a (\Omega _n \cap D'')) + \mu (\Omega_{n+1}) \;\; , \; n\geq 0.$$
If we denote $W_n =f^{-p}_a (\Omega _n \cap D'')$ then we get inductively that
$$\mu (\Omega _0) \geq \sum_{n\geq 0 } \mu (W_n ).$$
Since $|(f^p)'|_\sg$ is bounded on a $B(\{a\}\cup \post , 2T)$, there exists $L>1$ such that
$W_n \supset D(a,L^{-n} )$ for every $n\geq 1$. Therefore
$$\mu (\Omega _0) \geq \sum_{n\geq 0 } \mu ( D(a,L^{-n} )) \geq \sum_{n\geq 0 } \mu ( f^{-1} (D(a,L^{-n} )) \cap U_a)$$
with $U_a$ a logarithmic tract over the asymptotic value $a$. But on $U_a$ $\mu$ is equivalent
to the conformal measure $m$ (Lemma \ref{4.3}) and, with the same calculations that lead to (\ref{3.7}), we get that
\beq \label{4.8} \sum_{n\geq 0} m( f^{-1} (D(a,L^{-n} )) \cap U_a)
\asymp \sum_{n\geq 0} \left(\sum_{N\geq n }\sum_{k}  \left(N^2+k^2 \right)^{-\frac{\rho+1}{2\rho } h}\right) \eeq
which is finite if and only if $h> 3 \frac{\rho}{\rho +1}$.
\epf

\

It remains to investigate the case $h>3\frac{\rho}{\rho+1}$. In order to do so we write
\beq \label{4.9}
f^{-p} (\Gamma _n ) =\Gamma _{n+1} \cup W_n\cup S_n
\eeq
where
$$W_n= \bigcup_{a\in\av }W_n^a  \quad with \quad W_n^a = f_a^{-p} (D''_a \cap \Gamma _n ), \;\;  D''_a=D(f^p(a),T)$$
 and where $S_n$ is the remaining set.
The measure $\mu$ being $f$--invariant, the sequence $(\mu(\Gamma _n) )_n$ is decreasing.
We need the following additional property.

\

\blem \label{4.10}
For the $\sg$--finite invariant measure $\mu$ we have that $\lim _{n\to\infty } \mu(\Gamma _n) =0$.
\elem

\bpf Let $l=\lim _{n\to\infty } \mu(\Gamma _n) $.
From (\ref{4.9}) follows inductively that
$$\mu(\Gamma _0)= l + \sum_{n=0}^\infty (\mu (W_n ) +\mu(S_n)).$$
It is therefore natural to consider the set
$B=\bigcup_{n=0}^\infty (W_n\cup S_n)$. Define $\Gamma_\infty = \Gamma _ 0\cup B$ and
let $f_\infty $ be the induced map, i.e. the first return map, of $f^p $ on the set $\Gamma_\infty$.
Since $\mu$ is conservative, the conditional measure $\mu_\infty =\frac{\mu }{\mu (\Gamma _\infty )}$
is $f_\infty$ invariant (see \cite{aa}). Hence
$$ \mu_\infty (B) =\mu_\infty (f_\infty ^{-1} (\Gamma _0)) =\mu_\infty (\Gamma _0) $$
which implies that $\mu(B)=\mu(\Gamma _0)$. But this is only possible if $l=0$.
\epf

\

\ni The last step of the proof of Theorem \ref{4.5} is the following.

\

\blem\label{4.11}
If $h_f>3 \frac{\rho}{\rho +1}$ then the measure $\mu$ is finite.
\elem

\

\bpf
We have to show that $\mu (\Omega _0)< \infty$.
Since $\lim _{n\to\infty } \mu(\Gamma _n) =0$ it follows from induction that
$$
\mu (\Omega _0)
=\sum_{N=0}^\infty \mu ( \Gamma _N)
= \sum_{N=0}^\infty \left( \sum_{n=N}^\infty \mu(W_n)+\mu(S_n) \right).
$$
Let us first consider the term corresponding to $S_n$.

Choose again a Besicovitch covering of $\Omega_0$ by discs $D_j=D(x_j,2T)$, $x_j\in \post$. Let $D$ be one
of these discs and denote by $f_*^{-p}$ the inverse branches of $f^p$ defined on $D$ such that
$$S_n =\bigcup_{D\in \{D_j\}} \bigcup_* f_*^{-p} (\Gamma _n \cap D) \quad for \;\; every \;\; n\geq 0.$$
Since there is $c>0$ for which the sets
 $S_n \subset \Delta = \oc \setminus B(\av \cup \post , cT)$
we have $\mu(S_n ) \asymp m(S_n)$ (Lemma \ref{4.3}). Therefore we can do the following estimation.
$$\mu \left( \bigcup_* f_*^{-p} (\Gamma _n \cap D)\right) \asymp \sum_* m (f_*^{-p} (\Gamma _n \cap D))
\asymp \sum_* |(f^p)'(z_*)|_\sg ^{-h} m(\Gamma_n \cap D)$$
where, for every $*$, $z_*$ is any fixed point in $f_*^{-p} (D)$.
Since $D\cap B(\av , T)=\emptyset$ it follows from Lemma \ref{2.1} together with Proposition \ref{2.0.2}
that
$$\mu \left( \bigcup_* f_*^{-p} (\Gamma _n \cap D)\right)
\preceq \sum_* (1+|z|^{\rho+1}) ^{-h} m(\Gamma_n \cap D)
\preceq m(\Gamma_n \cap D).$$
Summing now over the discs of the Besicovitch covering and using the exponential decay of the
$m$-mass of the sets $\Gamma _n$ given in Lemma \ref{3.10} we finally get
$$
\mu (S_n) \preceq  m(\Gamma _n) \preceq \g ^n ,$$
and thus
$$ \sum_{N=0}^\infty \sum_{n\geq N}^\infty \mu (S_n) <\infty .$$

It suffices now to obtain the corresponding statements for the sets $W_n$.
Notice again that $\mu(W_n)= \mu ( f^{-1} (W_n ))\asymp m( f^{-1} (W_n ))$.
The set $ f^{-1} (W_n )$
contains a subset that lies in parabolic tracts and a remaining set
say $S_n'$. The $m$--mass of the later can be estimated exactly like we just did for $S_n$.
It therefore suffices to see what happens in just one tract $U_a$ and to estimate
the mass of $U_a \cap f^{-1} (W_n )$.
Clearly there is $c>0$ such that $W_n \subset D(a, c 2^{-n})$. We therefore can conclude
precisely like in (\ref{4.8}) that
$$\sum_{N=0}^\infty \sum_{n\geq N}^\infty \mu (W_n) <\infty $$
if and only if $h>3\frac{\rho}{\rho +1}$.
\epf

\section{Bowen's formula, Hausdorff dimension and Hausdorff measures}
We start with the following fact concerning the $h$-dimensional Hausdorff measure
$\cH^h$ on $\jul$.

\

\bprop\label{p1v83}
If $h<2$, then $h$-dimensional Hausdorff measure of $\jul$ vanishes, $\cH^h(\jul)=0$.
If $h=2$, then $\jul=\oc$. In either case $\HD(\jul)\le h$.
\eprop

\

\bpf Fix an arbitrary $z\in\jul\sms\bu_{n=0}^\infty f^{-n}(\av\cup\{\infty\})$. Then
there exists an increasing unbounded sequence $(n_j)_{j=1}^\infty$ such that for every
$j\ge 1$ there exists a meromorphic inverse branch $f_z^{-n_j}:D(f^{n_j}(z),2T)\to\oc$
sending $f^{n_j}(z)$ to $z$. Then $f_z^{-n_j}\(D(f^{n_j}(z),T)\)\sbt
D\(z,KT|(f^{n_j})'(z)|^{-1}\)$, and therefore
$$
\aligned
m\(D\(z,KT|(f^{n_j})'(z)|_\sg^{-1}\)\)
&\geq K^{-h} |(f^{n_j})'(z)|_\sg^{-h}m\(D(f^{n_j}(z),T)\)\\
&\geq M_m(T) (K^2T)^{-h} \(KT|(f^{n_j})'(z)|_\sg^{-1}\)^h.
\endaligned
$$
Hence,
$$
\limsup_{r\to 0}{m(D(z,r))\over r^h}
\ge \liminf_{j\to\infty}{m\(D\(z,KT|(f^{n_j})'(z)|_\sg^{-1}\)\)
     \over\(KT|(f^{n_j})'(z)|_\sg^{-1}\)^h}
\ge M_m(T)(K^2T)^{-h}>0.
$$
Thus
\beq\label{1v65}
\cH^h|_{\jul}\le Cm
\eeq
with some universal constant $C>0$. Proceeding further, suppose first that $h<2$. Recall that
$W_R=\{z\in\oc:|z|>R\}$. It follows from (\ref{3.4a}) that
\beq\label{2v83}
m(W_R)\gek R^{2\rho -(\rho +1)h}.
\eeq
Due to conservativity and ergodicity of the measure $m$, there exists a Borel set
$Y\sbt \jul\sms \bu_{n\ge 0}f^{-n}(\infty)$ such that $m(Y)=1$ and $\infty\in\om(z)$
for all $z\in Y$. Fix one $z\in Y$. There then exists an unbounded increasing sequence
$(n_j)_1^\infty$ such that
\beq\label{1v85}
\lim_{j\to\infty}|(f^{n_j})(z)|=+\infty \  \text{ and } \  |(f^{n_j})(z)|\ge 4T^{-1}
\eeq
for all $j\ge 1$. So there exist meromorphic inverse branches $f_z^{-n_j}:W_{|(f^{n_j}(z)|}
\to\oc$ sending $f^{n_j}(z)$ to $z$. Put $r_j=2K|(f^{n_j})'(z)|_\sg^{-1}|f^{n_j}(z)|^{-1} $.
Looking at (\ref{2v83}), we get
$$
\aligned
{m(D(z,r_j))\over r_j^h}
&\ge {m\(f_z^{-n_j}\(D(f^{n_j})(z),2|(f^{n_j})(z)|^{-1})\)\)\over r_j^h} \\
&\ge {K^{-h}|(f^{n_j})'(z)|_\sg^{-h}m\(D(f^{n_j})(z),2|(f^{n_j})(z)|^{-1})\)\over
      K^h   |(f^{n_j})'(z)|_\sg^{-h}|(f^{n_j})(z)|^{-h}} \\
&\gek K^{-2h}|(f^{n_j})(z)|^hm\(W_{|(f^{n_j})(z)|}\)
 \ge  K^{-2h}|(f^{n_j})(z)|^h|(f^{n_j})(z)|^{2\rho -(\rho +1)h} \\
&=    K^{-2h}|(f^{n_j})(z)|^{\rho (2-h)}.
\endaligned
$$
Since $2-h>0$, we therefore conclude from this and (\ref{1v85}) that
$$
\limsup_{r\to 0}{m(D(z,r))\over r^h}
\ge \lim_{j\to\infty}{m(D(z,r_j))\over r_j^h}
\ge \lim_{j\to\infty}K^{-2h}|(f^{n_j})(z)|^{\rho (2-h)}
=+\infty.
$$
Thus, $\cH^h(Y)=0$. Since $\cH^h(\jul\sms Y)=0$ by (\ref{1v65}), we thus have
$\cH^h(\jul)=0$. The case $h<2$ is done. If $h=2$, then for the sequence
$(n_j)_1^\infty$ from the beginning of the proof, we will have $m(D(z,r_j))\comp
r_j^2$, which implies that $m$ and $l_s$, the spherical Lebesgue measure on
$\oc$, are equivalent. So, $l_s(\jul)>0$. Now, if $\jul\ne\oc$, then $\jul$ would
be nowhere dense in $\oc$, and in the same way as Lemma~\ref{l3v54}, making
use of the Lebesgue Density Theorem, we would prove that $l_s(\jul)=0$. This
contradiction finishes the proof.
\epf

\

\ni Although $\cH^h(\jul)=0$ (if $h<2$), we shall however show that $h=\HD(\jul)$.
The proof will utilize the induced (first return) map we are going to describe now.
Let
\beq\label{1v69}
X=\jul\sms\(B(\post,T)\cup\bu_{a\in\av}f_a^{-1}(D(f(a),T))\).
\eeq
Let $f_*:X\to X$ be the first return map of $f$ on $X$. That is
$$
f_*(x)=f^{\tau(x)}(x),
$$
where $\tau(x)=\min\{n\ge 1:f^n(x)\in X\}$. Since $f:\jul\to\jul$ is conservative
with respect to the measure $\mu$ (see Theorem~\ref{t1v55}), the map $f_*$ is
well-defined on the complement of a set of $\mu$ measure zere, in fact, as it is
easy to see, it is well-defined on the complement of $\bu_{n\ge 0}f^{-n}(\post)$,
which is of measure zero by Lemma~\ref{l3v54} and by formula (\ref{4.2}). Since
the Radon-Nikodym derivative ${d\mu/dm}$ is uniformly bounded above pn $X$, $\mu(X)
<+\infty$. For every $x\in X$ define
$$
f_*'(x)=\(f^{\tau(x)}\)'(x) \  \text{ and }  \
|f_*'(x)|_\sg=\lt|\(f^{\tau(x)}\)'(x)\rt|_\sg.
$$
We shall prove the following.

\

\blem\label{l1v71}
$\b:=\inf\{|f_*'(z)|_\sg:z\in X\}>0$ and there exists $k\ge 1$ so large that
$|(f_*^k)'(z)|_\sg\ge 2$ for all $z\in X$.
\elem

\bpf
In the course of the proof of this lemma $Q$ stands for an appropriately large positive
constant.

Suppose first that $z\in X\cap U_a$ where $U_a$ is a logarithmic tract over some $a\in \av$
such that $f(U_a)=f_a^{-1} (D(f(a),T))$. Let $n\ge 0$ the least integer such that
$f^{n+1}(z)\notin D(\post,T)$. Then
\beq\label{3v71}
\aligned
|f_*'(z)|_\sg
&=|(f^{n+1})'(z)|_\sg
\ge Q^{-1}|f(z)-a|(1+|z|^{\rho+1})|(f^n)'(z)|_\sg \\
&\ge Q^{-2}(1+|z|^{\rho+1})
\ge Q^{-2}
\endaligned
\eeq
and
\beq\label{4v71}
|f_*'(z)|_\sg \ge 2Q^2
\quad
\text{if, in addition,} \;\; |z|\ge R.
\eeq
 For all other $z\in X$,
Lemma~\ref{2.1} implies that
\beq\label{1v71}
|f_*'(z)|_\sg\ge Q^{-1}\ge Q^{-2}.
\eeq
If in addition $|z|>R$ with $R>0$ large enough, then
\beq\label{2v71}
|f_*'(z)|_\sg\ge 2Q^2.
\eeq
 The first part of our lemma is thus proved. We shall now
demonstrate the following.

\sp\ni Claim: There exists $l=l(R)\ge 1$ such that
$$
|(f_*^n)'(z)|_\sg\ge 2Q^2
$$
for all $n\ge l$ and all $z\in D(0,R)\cap X$.

\ni Indeed, suppose for the contrary that there exist an increasing sequence $n_j\to\infty$
and a sequence $z_j\in X\cap\ov D(0,R)$ such that
\beq\label{1v73}
|(f_*^{n_j})'(z_j)|_\sg<2Q^2
\eeq
for all $j\ge 1$. Since $f_*^{n_j}(z_j)\in X$, there exists a unique meromorphic
inverse branch $f_{z_j}^{-N_j}:D(f^{n_j}(z_j,2T)\to\oc$ of $f^{N_j}$, sending
$f_*^{n_j}(z_j)$ to $z_j$, where $N_j=\tau(z_j)+\tau(f_*(z_j))+\ld+\tau(f_*^{n_j-1}(z_j))$.
It then follows from (\ref{bdp}) and (\ref{1v73}) that
$$
f_{z_j}^{-N_j}\(D(f_*^{n_j}(z_j,T/2))\)\sbt D(z_j,(4KTQ^2)^{-1}),
$$
or equivalently,
$$
f^{N_j}\(D(z_j,(4KTQ^2)^{-1})\)\sbt D\(f_*^{n_j}(z_j),T/2\).
$$
Passing to a subsequence we may assume without loss of generality that the
sequnce $(z_j)_1^\infty$ converges to a point $z\in\jul\cap\ov D(0,R)$ and
$|z_j-z|<(8KTQ^2)^{-1}$ for all $j\ge 1$. Since $D(f_*^{n_j}(z_j,T/2)\cap
B(\post,T/2)=\es$, it follows from Montel's Theorem that the family
$\lt\{f^{N_j}|_{D(z,(8KTQ^2)^{-1})}\rt\}_{j=1}^\infty$ is normal, contrary to
the fact that $z\in\jul$. The claim is proved.

\sp\ni Let $k=2l$. If $|f_*^j(z)|\ge R$ for all $j=1,2,\ld,l$, then by (\ref{1v71}) -
(\ref{4v71}), we get
$$
|(f_*^k)'(z)|_\sg\ge (2Q^2)^lQ^{-2l}=2^l\ge 2.
$$
If $|f_*^j(z)|< R$ for some $0\le j\le l$, then let $j$ be minimal with this
property. It then follows from (\ref{2v71}), (\ref{4v71}) and the claim, that
$$
|(f_*^k)'(z)|_\sg
=|(f_*^j)'(z)|_\sg|(f_*^{k-j})'(f^j(z))|_\sg
\ge |(f_*^{k-j})'(f^j(z))|_\sg
\ge 2Q^2
\ge 2.
$$
We are done.
\epf

\

\ni Now, we shall prove the following.

\

\blem\label{l1v73}
The function $z\mapsto\log|f_*'(z)|_\sg$ is integrable on $X$ with respect to $\mu_X$,
the conditional measure on $X$ induced by $\mu$. In addition
$\chi:=\int\log|f_*'|_\sg d\mu_X>0$.
\elem

\bpf
Since the Radon Nikodym derivative $d\mu/dm$ is uniformly bounded on $X$, it
suffices to demonstrate that the function $z\mapsto\log|f_*'(z)|_\sg$ is integrable
on $X$ with respect to the measure $m$ ($\chi>0$ follows immediately from
Lemma~\ref{l1v71}). For every $a\in\av$ let $A_n(a)$, $n\ge 0$, be the annuli
defined by formula (\ref{3.5}). Put $A_n=\bu_{a\in\av}A_n(a)$. Partition
$X\sms f(A_0)$ by disjoint Borel sets $X_n$, $n\ge 0$, such that
$D(x_n,2\diam(X_n))\cap (\av\cup\post)=\es$ with some $x_n\in X_n$. Then, because of
Lemma~\ref{2.1} and Proposition~\ref{2.0.2}, we get that
\beq\label{1v75}
\aligned
\int_{X\cap f^{-1}(X\sms f(A_0))} &\lt|\log|f_*'|_\sg \rt| dm = \\
&    =\sum_{n=1}^\infty\int_{X\cap f^{-1}(X_n)}\lt| \log|f_*'|_\sg \rt| dm \\
&\comp\sum_{n=1}^\infty m(X_n)\sum_{z\in X\cap f^{-1}(w_n)}|f_*'(z)|_\sg^{-h}\lt|\log|f_*'(z)|_\sg\rt| \\
&\comp\sum_{n=1}^\infty m(X_n)\sum_{z\in X\cap f^{-1}(w_n)}(1+|z|^{\rho+1})^{-h}
      \lt|\log(1+|z|^{\rho+1})+O(1)\rt| \\
&\lek \sum_{n=1}^\infty m(X_n)\sum_{z\in X\cap f^{-1}(w_n)}(1+|z|^{\rho+1})^{-t} \\
&\le M_t\sum_{n=1}^\infty m(X_n)
 \le M_t
 <+\infty,
\endaligned
\eeq
where $w_n$ is an arbitrary point in $X_n$ and $t$ is a fixed number in
$({\rho\over \rho+1},h)$. Now, following notation from Prposition~\ref{3.4}, for every
$a\in\av$ and every $n\ge 0$, set
$$
\Ga_a=f^{-1}(f(a))\sms (\av\cup\post),
$$
$$
Y_n(a)=\bu_{b\in\Ga}f_b^{-1}\circ g_{f(a)}^{-n}(A_n(a))\cup
       \bu_{b\in f^{-1}(a)}f_b^{-1}\circ f_a^{-1}\circ g_{f(a)}^{-n}(A_n(a))
$$
and
$$
Y_a=\bu_{n=0}^\infty Y_n(a).
$$
Keep $t\in({\rho\over \rho+1},h)$. Again, in virtue of Lemma~\ref{2.1} and
Proposition~\ref{2.0.2}, and also of Lemma~\ref{3.10}, we get that
$$
\aligned
\int_{Y_a} & \lt|\log|f_*'|_\sg\rt| dm =\\
&    =\sum_{n=0}^\infty\int_{Y_n(a)}\lt|\log|f_*'|_\sg\rt| dm \\
&\comp\sum_{n=0}^\infty\bigg(\sum_{b\in\Ga}m\(f_b^{-1}\circ g_{f(a)}^{-n}(A_n(a))\)
                       \lt|\log|f'(b)||(g^n)'(f(a))|+O(1)\rt| + \\
& \  \  \  \  \     +\sum_{b\in f^{-1}(a)}m\(f_b^{-1}\circ f_a^{-1}\circ g_{f(a)}^{-n}(A_n(a))\)
                       \lt|\log|f'(b)||(g^n)'(f(a))|+O(1)\rt|\bigg) \\
&\lek \sum_{n=0}^\infty\sum_{b\in\Ga}(1+|b|^{\rho+1})^{-h}\g^n\lt|\log(1+|b|^{\rho+1})
        +\log|(g^n)'(f(a))|+O(1)\rt| + \\
& \  \  \  \  \    +\sum_{n=0}^\infty\sum_{b\in f^{-1}(a)}(1+|b|^{\rho+1})^{-h}\g^n
       \lt|\log(1+|b|^{\rho+1})+\log|(g^n)'(f(a))|+O(1)\rt| \\
&\le  \sum_{n=0}^\infty\g^n\sum_{b\in\Ga}(1+|b|^{\rho+1})^{-h}\g^n\lt|\log(1+|b|^{\rho+1})+O(n)\rt| + \\
& \  \  \  \  \     +\sum_{n=0}^\infty\g^n\sum_{b\in f^{-1}(a)}
          (1+|b|^{\rho+1})^{-h}\g^n\lt|\log(1+|b|^{\rho+1})+O(n)\rt|\\
&\lek  \sum_{n=0}^\infty\g^n\sum_{b\in\Ga\cup f^{-1}(a)}(1+|b|^{\rho+1})^{-h}\lt|\log(1+|b|^{\rho+1})\rt|
      +\sum_{n=0}^\infty n\g^n\sum_{b\in\Ga\cup f^{-1}(a)}(1+|b|^{\rho+1})^{-h} \\
&\lek  \sum_{n=0}^\infty\g^n\sum_{b\in\Ga\cup f^{-1}(a)}(1+|b|^{\rho+1})^{-t}
      +M_h\sum_{n=0}^\infty n\g^n \\
&<+\infty.
\endaligned
$$
Hence,
\beq\label{1v77}
\int_{\bu_{a\in\av}Y_a}\lt|\log|f_*'|_\sg\rt| dm <+\infty.
\eeq
Finally, for every $a\in\av$, let
$$
U_a=\bu_{n,k\ge 1}U_{n,k}(a).
$$
In view of (\ref{3.4a}) and (\ref{3.4b}) we get that
$$
\aligned
\int_{U_a}\lt|\log|f_*'|_\sg\rt| dm
&     =\sum_{n,k\ge 1}\int_{U_{n,k}(a)}\lt|\log|f_*'|_\sg\rt| dm \\
&\comp\sum_{n,k\ge 1}m(U_{n,k}(a))\lt|\log\(|z_{n,k}|^{\rho+1}|f(z_{n,k})-a||f(z_{n,k})-a|^{-1}\)\rt| \\
&\comp\sum_{n,k\ge 1}|z_{n,k}|^{-h(\rho+1)}\lt|\log\(|z_{n,k}|^{\rho+1}\rt| \\
&\comp\sum_{n,k\ge 1}(n^2+k^2)^{-h{\rho+1\over 2\rho}}\log(n^2+k^2) \\
&\lek  \sum_{n,k\ge 1}(n^2+k^2)^{-t{\rho+1\over 2\rho}}\\
&< +\infty.
\endaligned
$$
Hence, $\int_{\bu_{a\in\av} U_a}\lt|\log|f_*'|_\sg\rt| dm<+\infty$. Summing up this,
(\ref{1v75}), and (\ref{1v77}), we conclude that $\int_X\lt|\log|f_*'|_\sg\rt| dm<+\infty$,
and the proof is complete.
\epf

\

\ni The main result of this section is this.

\

\bthm\label{t1v81}
It holds $\HD(\jul)=h$.
\ethm

\bpf
In view of Proposition~\ref{p1v83} it suffices to show that $\HD(\jul)\ge h$. Let
$X\sbt\jul$ be the set defined by (\ref{1v69}) and let $f_*:X\to X$ be the corresponding
induced map. In virtue of Lemma~\ref{l1v73}, Lemma~\ref{l1v71} and Birkhoff's Ergodic
Theorem, there exists a Borel set  $\hat X\sbt X$ such that $\mu(\hat X)=1$ and
$$
\lim_{n\to\infty}{1\over n}\log|(f_*^n)'(z)|_\sg=\chi>0
$$
for all $z\in\hat X$. In particular,
\beq\label{1v81}
\lim_{n\to\infty}{\log|(f_*^{k(n+1)})'(z)|_\sg\over\log|(f_*^{kn})'(z)|_\sg}=1,
\eeq
where $k\ge 1$ comes from Lemma~\ref{l1v71}. For every $z\in\hat X$ and every $n\ge 0$
define
$$
r_n(z)=(2K)^{-1}|(f_*^{kn})'(z)|_\sg^{-1}.
$$
Fix $\e\in(0,h)$. In virtue of (\ref{1v81}), for every $z\in\hat X$ we have
\beq\label{1v80}
{r_n(z)\over r_{n+1}(z)}\le r_n(z)^{-{\e\over 2}}
\eeq
for all $n\ge 1$ large enough. It follows from (\ref{bdp}) and conformality of $m$
that
\beq\label{1v83}
\aligned
m(D(z,r_n))
&\le m\(f_*^{-N_{kn}(z)}\(D(f_*^{kn}(z),T/2)\)\) \\
&\le K^h\lt|\lt(f_*^{N_{kn}(z)}\rt)'(z)\rt|_\sg^{-h}m\(D(f_*^{kn}(z),T/2)\) \\
&\le K^h\lt|\lt(f_*^{kn}\rt)'(z)\rt|_\sg^{-h}
 =   (2K^2T)^hr_n^h.
\endaligned
\eeq
Now, keeping $z\in\hat X$, take an arbitrary radius $r\in(0,(2K)^{-1}T)$. Since the
sequence $(r_n)_0^\infty$ is strictly decreasing, there exists a unique $n\ge 0$
such that $r_{n+1}\le r<r_n$. In view of (\ref{1v83}) and (\ref{1v80}), we get that
\beq
\aligned
\lim_{r\to 0}{m(D(z,r))\over r^{h-\e}}
&\le\lim_{n\to\infty}{m(D(z,r_n))\over r_{n+1}^{h-\e}}
 =  \lim_{n\to\infty}\lt({m(D(z,r_n))\over r_n^{h-\e}}\lt({r_n\over r_{n+1}}\rt)^{h-\e}\rt) \\
&\le\lim_{n\to\infty}\(r_n^\e r_n^{-{\e\over 2}}\)
 =  \lim_{n\to\infty}r_n^{{\e\over 2}}
 =  0.
\endaligned
\eeq
Since $m(\hat X)>0$, we therefore conclude that $\cH^{h-\e}(\jul)\ge \cH^{h-\e}(\hat X)
=+\infty$. Thus $\HD(\jul)\ge h-\e$, and eventually, letting $\e\downto 0$, we get
$\HD(\jul)\ge h$. We are done.
\epf

%********************************************************************************************************

\end{document}